\documentclass[11pt]{article}

\usepackage[utf8]{inputenc}
\usepackage[english]{babel}
\usepackage{amsmath,amsthm,amsfonts}
\usepackage{graphicx}
\usepackage{bm}
\usepackage{tikz}
\usepackage{multirow}
\usepackage[makeroom]{cancel}
\usepackage{url}

\newtheorem{remark}{Remark}

\begin{document}
\title{Generalized Multiscale Inversion for Heterogeneous Problems}
\author{
Eric T. Chung\thanks{Department of Mathematics, The Chinese University of Hong Kong, Shatin, New Territories, Hong Kong SAR, China (\texttt{tschung@math.cuhk.edu.hk}) }
\and
Yalchin Efendiev\thanks{Department of Mathematics \& Institute for Scientific Computation (ISC),
Texas A\&M University,
College Station, Texas, USA (\texttt{efendiev@math.tamu.edu})}
\and
Bangti Jin\thanks{Department of Computer Science, University of College London, Gower Street, London WC1E 6BT, UK (\texttt{bangti.jin@gmail.com})}
\and
Wing Tat Leung\thanks{Department of Mathematics, Texas A\&M University, College Station, TX 77843, USA (\texttt{leungwt@math.tamu.edu})}
\and
 Maria Vasilyeva\thanks{Institute for Scientific Computation, Texas A\&M University, College Station, TX, USA \& Department of Computational Technologies, North-Eastern Federal University, Yakutsk, Republic of Sakha (Yakutia), Russia
(\texttt{vasilyevadotmdotv@gmail.com})}
}

\newcommand{\rlg}[1]{}

\maketitle
\begin{abstract}
In this work, we propose a generalized multiscale inversion algorithm for heterogeneous problems that aims at solving
an inverse problem on a computational coarse grid. Previous inversion techniques for multiscale problems seek a coarse-grid
media properties, e.g., permeability and conductivity, and by doing so, they assume that there exists a homogenized
representation of the underlying fine-scale permeability field on a coarse grid. Generally such assumptions do not hold for highly
heterogeneous fields, e.g., fracture media or channelized fields, where the width of channels are very small compared to
the coarse-grid sizes. In these cases, grid refinement can lead to many degrees of freedom, and thus unattractive to apply. The proposed algorithm is based on
the Generalized Multiscale Finite Element Method (GMsFEM), which uses local spectral problems to identify non-localized
features, i.e., channels (high-conductivity inclusions that connect the boundaries of the coarse-grid block). The inclusion
of these features in the coarse space enables one to achieve a good accuracy. The approach is valid under the assumption
that the solution can be well represented in a reduced-dimensional space by multiscale basis functions. In practice, these
basis functions are non-obervable as we do not identify the fine-scale features of the permeability field. Our inversion
algorithm finds the discretization parameters of the resulting system. By doing so, we identify the appropriate coarse-grid
parameters representing the permeability field instead of fine-grid permeability field.
We illustrate the approach by numerical results for fractured media.
\end{abstract}

\section{Introduction}

In many applications, one deals with medium properties of multiple scales and high contrast. For example, in subsurface
applications, high-conductivity channels or fractures can appear in multiple locations and have complex geometries. Such
features typically have multiple scales, e.g., very small widths and multiple (long) length scales. The related inverse
problems include finding permeability (or channel distribution) from noisy and sparse pressure or concentration
measurements, and they can be posed as a regularized least squares formulation and/or within a Bayesian formulation.

There are several challenges when performing inversion using standard approaches (see the monographs \cite{EnglHankeNeubauer:1996,
Tarantola:2005,KaipioSomersalo:2005,SchusterKaltenbacher:2012,ItoJin:2015} for a few references)
for heterogeneous problems. Because of the presence of small scales, one needs to resolve multiple scales properly, which
can lead to huge ill-posed systems that are difficult to solve. However, one cannot perform inversion on a coarse grid
using standard approaches directly, since the latter implicitly assumes that there is a homogenized model (see
e.g., \cite{NolenPavliotisStuart:2012,GullikssonHolmbom:2016,FrederickEngquist:2017} for related inverse problems for homogenization). It was shown
in \cite{GalvisEfendiev:2010a,EfendievGalvisHou:2013,ChungEfendievHou:2016} that this assumption is not valid for many
practical multiscale problems, even at a low-order approximation. Indeed, because of the presence of high-contrast channels,
one cannot use a single permeability or conductivity to represent a coarse-grid block. To remedy these drawbacks, multiple
continuum approaches \cite{ArbogastDouglasHornung:1990,BarenblattZheltovKochina:1960,KazemiMerrillPorterfieldZeman:1976,
PruessNarasimhan:1982,WarrenRoot:1963,WuPruess:1988} can be used in this context; however, these approaches require
multiple assumptions \cite{ChungEfendievLeungVasilyeva:2017}. Meanwhile, using fine-grid discretizations can lead to
many degrees of freedom without a priori knowledge of the locations of these thin features. In this paper, we present a novel
generalized multiscale inversion algorithm, which employs our recent multiscale methods and solves inverse problem for
discretization parameters rather than for fine-grid permeability fields. Thus by construction, it provides a
low-dimensional inverse problem on the coarse grid and avoids many prior assumptions on the fine-grid geometry in order to regularize
the inverse problem (in the spirit of regularization by discretization).

\begin{figure}[ht!]
\centering
\includegraphics[width=3in]{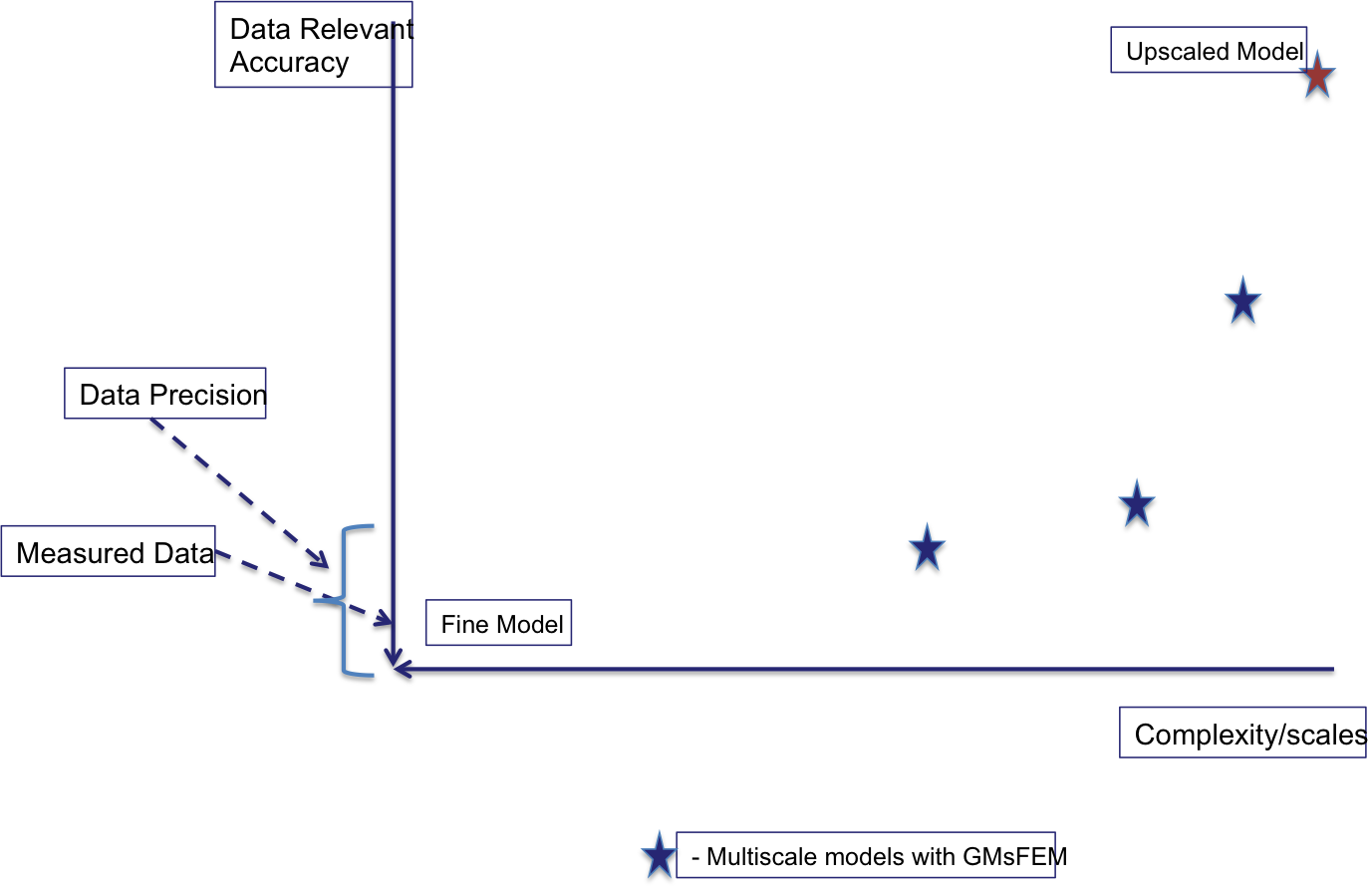}
\caption{A schematic illustration of the concept of multiscale inversion: The plot shows that one needs appropriate coarse-grid models (with multiple
basis functions) for the inversion in order to achieve an accuracy within the data error.}
\label{fig:chart}
\end{figure}

Next, we briefly discuss generalized multiscale methods in the context of inverse problems. We conceptually sketch it in
Fig. \ref{fig:chart}, where we emphasize that one needs appropriate coarse-grid models (with multiple basis functions)
for the inversion in order to achieve an accuracy within the error tolerance of the data. For simplicity, we consider a
multiscale parabolic equation
\begin{equation}\label{eq:original}
{\frac{\partial u}{\partial t}} - \text{div} (\kappa\nabla u) =f,\quad \mbox{in }\Omega\times(0,T],
\end{equation}
with a homogeneous Dirichlet boundary condition and a suitable initial condition, where $\Omega\subset\mathbb{R}^n$ is an
open bounded domain, $T>0$ is a fixed a final time, and $\kappa_0\leq\kappa\leq \kappa_1$ is the unknown permeability
field that varies over multiple scales with high contrast. Our approach begins with a computational grid, called the
coarse grid, which, as usual, does not resolve all the features of the permeability $\kappa(x)$. One standard
approach is to seek $\kappa^*(x)$ on a coarse grid directly. However, it automatically assumes that one has a homogenization
within a set of permeability fields that we seek. The latter assumption is violated in many important practical applications,
including, e.g., identifying fractures (thin high-conductivity features) or channels with extremely low or high conductivities.
In these cases, when the thin features are subgrid with respect to the coarse-grid block, homogenization can only
provide very inaccurate solutions. Some alternative approaches include multi-continuum, where multiple homogenized coefficients
are assigned in each block, which, however, need certain modeling assumptions. In this work, we shall employ generalized
multiscale approaches, where one constructs multiple physically-relevant basis functions in each coarse block from the
observational data (in an adaptive manner).

The multiscale method that we employ for the inversion is based on the Generalized Multiscale Finite Element Method (GMsFEM)
\cite{HouWu:1997,EfendievGalvisWu:2011,EfendievGalvisHou:2013,GalvisEfendiev:2010a,ChungEfendievHou:2016}. The main idea
of the GMsFEM is to construct multiscale basis functions in each coarse block, by solving local spectral problems. The
multiscale basis functions are selected based on dominant modes of local spectral problems. The dominant modes can be
identified through a spectral gap and the dominant modes correspond to channelized features, i.e., the high-conductivity
channels that connect the boundaries of the coarse block. These features cannot be localized and require separate basis
functions. If these features are not represented by separate basis functions and represented by fewer basis functions,
one can only get very inaccurate solutions. Hence, if one uses only an upscaled permeability (which corresponds to one
basis function), the inversion can provide an inaccurate solution.

Our generalized multiscale inversion algorithm formulates the inverse problem for the discretization parameters on a coarse
grid directly. The solution to the direct problem is assumed to be represented/captured by several basis functions in each
coarse block, where basis functions are not known \textit{a priori}, but to be inferred from the observational data simultaneously.
Next, we represent the measurements in terms of coarse-grid parameters, e.g., entries of the stiffness and mass matrices. The
latter is feasible under certain assumptions on physical nature of measurements. For example, if the measured quantities can be
written on a coarse grid, one can easily represent the observed data via coarse-grid parameters. Note that in our inversion
algorithm, we do not identify detailed basis functions, but only some average information that these basis functions will provide.
We call these multiscale basis functions {\it unobservable} and introduce {\it observable} counterpart, which allows extracting
some average information about the solution. Naturally, in the proposed algorithm, one needs certain physical constraints
(on the permeabilities etc.) in order to be able to recover some elements of stiffness and mass matrices. The proposed method
can also be formulated in a Bayesian framework, by imposing a prior
on the stiffness and mass matrices generated from a known fine-grid permeability field, and then to sample the resulting
posterior distribution with Markov chain Monte Carlo in order to quantify the associated uncertainties \cite{EfendievJinPreshoTan:2015}.

In the paper, we present several numerical examples for flows in fractured media, using a setup for shale gas
applications \cite{AkkutluEfendievVasilyeva:2016}, where the true model has fracture distributions that differ from
the initial model and the data are coarse-grid pressures. Because of fracture networks, we assume that the model
has at most two basis functions in each coarse block and perform inversion. We test the sensitivity of our
approach with respect to data noise and measurement location. Moreover, we present adaptive approaches,
where multiscale basis functions are used only in selected regions for the purpose of updating.

The rest of the paper is organized as follows. In Section \ref{sec:prelim}, we give some preliminaries about
grids, multiscale method, and the setup of the inverse problem. In Section \ref{sec:method}, we present our
generalized inversion algorithm. Numerical results are presented in Section \ref{sec:num}.

\section{Preliminaries}
\label{sec:prelim}
In this section, we describe preliminaries about generalized multiscale finite element methods (GMsFEM), and the setup
for the inverse problem.

\subsection{Coarse and fine grids}

\begin{figure}[ht!]
\centering
\includegraphics[width=0.6 \textwidth]{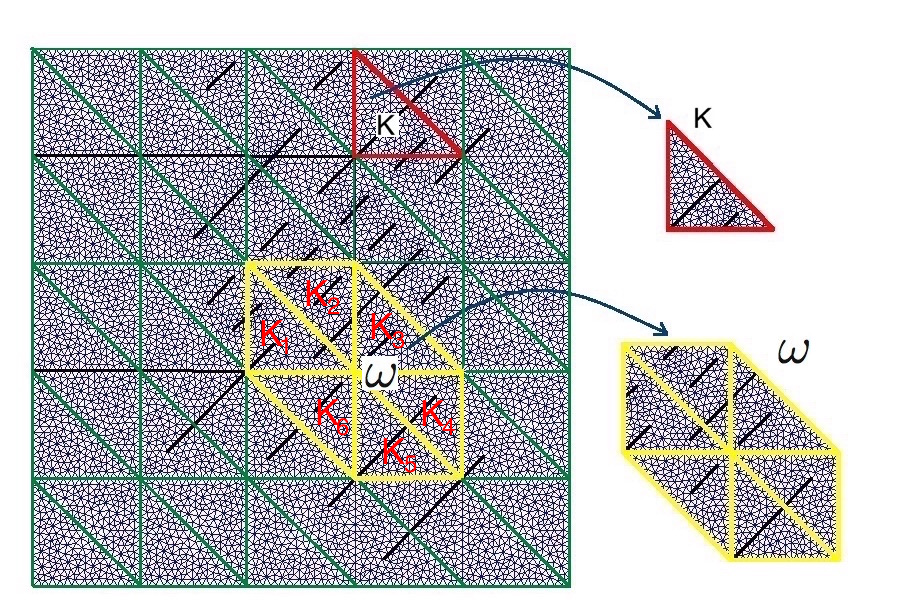}
\caption{Illustration of the coarse grid $\mathcal{T}^{H}$, coarse cell $K$, domain $\omega$ (the union of a few
coarse cells) and  fine grid $\mathcal{T}^{h}$.}
\label{fig:illustration}
\end{figure}

First we introduce the notion of fine- and coarse-grids. Let $\mathcal{T}^{H}$ be a conforming partition of the domain
$\Omega$ into finite elements, called coarse grid, with $H$ being the coarse-mesh size. Let $N_c$ be the number of
vertices, and $N$ the number of elements in the coarse mesh. Then each coarse element is further partitioned into a connected
union of fine-grid blocks, denoted by $\mathcal{T}^{h}$. The partition $\mathcal{T}^{h}$ is a refinement of the coarse grid
$\mathcal{T}^{H}$ with the mesh size $h$. Throughout, it is always assumed that the fine grid is sufficiently fine to resolve
the solution. We refer to Fig. \ref{fig:illustration} for an illustration.


\subsection{Multiscale basis functions}

The GMsFEM consists of two stages: offline and online. First we describe the online stage. Let $V = H^1_0(\Omega)$. Then the
solution $u$ of problem \eqref{eq:original} satisfies
\begin{equation}\label{eq:finesol}
({\frac{\partial u}{\partial t}},v) + a(u,v) = ( f,v) \quad \text{ for all } v\in V,
\end{equation}
where $a(u,v)=\int_{\Omega}\kappa \nabla u \cdot \nabla v\mathrm{d}x$, and $(\cdot,\cdot)$ denotes the $L^2$-inner product on $\Omega$.
Let $V_{ms}\subset V$ be the space spanned by all multiscale basis functions, whose construction is to be described in detail below.
Then the multiscale solution $u_{ms}$ is defined as: find $u_{ms}\in V_{ms}$ such that
\begin{equation}\label{eq:mssol}
a(u_{ms},v) = ( f,v) \quad \text{ for all } v\in V_{ms}.
\end{equation}

Next we describe the construction of the multiscale basis functions. It is performed on the fine mesh, even
though it is not use in our inversion. In the offline stage, a small dimensional finite element
space is constructed to solve the global problem for any input parameter, e.g., right-hand side or boundary
condition, on a coarse grid. The snapshot space $V_{H,\text{snap}}^{(i)}$ is constructed for a generic domain
$\omega_i$. The snapshot solutions are then used to compute multiscale basis functions. The ideal snapshot space
should provide a fast convergence and problem-relevant restrictions on the coarse spaces (e.g., divergence free
solutions), while can reduce the cost associated with constructing the offline spaces. One can generate
snapshot spaces in several different ways \cite{ChungEfendievHou:2016}, and here we employ harmonic snapshots
in an oversampling domain (cf. Fig. \ref{fig:illustration} for a sketch).

The snapshot space $V_{H,\text{snap}}^{(i)}$ consists of harmonic extensions of fine-grid functions that are defined on the boundary
$\partial\omega_i$. For each fine-grid function $\delta_l^h(x)$, we define $\delta_l^h (x_k)=\delta_{l,k},\,\forall
x_k\in \textsl{J}_{h}(\omega_i)$ ($\delta_{l,k}$ is the Kronecker symbol, i.e., $\delta_{l,k} = 1$ if $l=k$ and
$\delta_{l,k}=0$ if $l\ne k$), where the notation $\textsl{J}_{h}(\omega_i)$ denotes the set of fine-grid boundary nodes on
$\partial\omega_i$. Then we obtain a snapshot function $\eta_l^{(i)}$ by
\[
\mathcal{L}( \eta_{l}^{(i)})=0\ \ \text{in} \ \omega_i,\quad  \eta_{l}^{(i)}=\delta_l^h(x) \mbox{ on }\partial\omega_i.
\]
The snapshot functions can be computed in the oversampling region $\omega_i^{+}$ in order to enhance the convergence rate,
and one can use randomized boundary conditions to further reduce the associated cost \cite{CaloEfendievGalvisLi:2016},
in the spirit of randomized singular value decomposition.

The offline space $V_{ms}^{(i)}$ is computed for each $\omega_i$ (with elements of the space denoted $\psi_l^{(i)}$) from
the snapshot space $V_{H,\text{snap}}^{(i)}$. Specifically, we perform a spectral decomposition in the snapshot space and select the
dominant modes (corresponding to the smallest eigenvalues) to construct the offline (multiscale) space $V_{ms}^{(i)}$. The convergence
rate of the resulting method is determined by $1/\Lambda_*$, where $\Lambda_*$ is the smallest eigenvalue that the
corresponding eigenvector is not included in the multiscale space $V_{ms}^{(i)}$ \cite{EfendievGalvisWu:2011,Li:2017}.
The concrete formulation of the local spectral problem can be motivated from the error analysis as follows. The global energy
error can be decomposed into coarse subdomains. With the energy functional on the domain $\omega$ denoted by $a_\omega(u,u)$,
i.e., $a_\omega(u,u) =\int_\omega \kappa \nabla u\cdot \nabla u\mathrm{d}x$, we have
\begin{equation}
a_\Omega(u-u_H,u-u_H)\preceq
\sum_\omega a_\omega(u^\omega-u_H^\omega,u^\omega-u_H^\omega),
\end{equation}
where $\omega$ are coarse regions ($\omega_i$), and $u^\omega$ is the localization of the solution. The local spectral problem
is chosen to bound the local error $a_\omega(u^\omega-u_H^\omega,u^\omega-u_H^\omega)$. Ideally, we look for the subspace $V_{ms}^{\omega}$ such that
for any $\eta\in V_{H,\text{snap}}^{\omega}$, there exists a function $\eta_0\in V_{ms}^{\omega}$ such that
\begin{equation}\label{eq:off1}
a_{\omega}(\eta-\eta_0,\eta-\eta_0)\preceq {\delta}s_{\omega}(\eta-\eta_0,\eta-\eta_0),
\end{equation}
where $s_{\omega}(\cdot,\cdot)$ is an auxiliary bilinear form, which has to be chosen properly to ensure the desired approximation property {\color{red}\cite{}}. The main
empirical observation is that with the snapshot spaces chosen suitably, the smallest eigenvalues correspond to the channelized features
\cite{EfendievGalvisWu:2011,EfendievGalvis:2011}, and thus it enables our multiscale inversion technique.

\subsection{Setup of inverse problem}

In the paper, our goal is to find some average information about the solution $u_h(x)$ and the permeability field $\kappa(x)$ given
measured data, denoted by $d$. Since our multiscale inversion technique does not identify $\kappa(x)$ and the solution $u_h(x)$ directly,
we denote the integrated responses by $\kappa_{ms}(x)$ and $u_{ms}(x)$. In a Bayesian framework, we write the inverse problem as
\begin{equation*}
P(\kappa_{ms}(x), u_{ms}(x)|d)\propto
P(d|\kappa_{ms}(x), u_{ms}(x)) \pi(\kappa_{ms}(x))\pi(u_{ms}(x)),
\end{equation*}
where $P(d|\kappa_{ms}(x), u_{ms}(x))$ is the likelihood function, $\pi(\kappa_{ms}(x))$ is the prior on multiscale
discretization parameters related to the coarse-grid $\mathcal{T}^H$, and $\pi(u_{ms}(x))$ is the prior on the coarse-grid solution.
We will describe the likelihood function and these priors more precisely later on. For the data $d$, we will assume that we measure average pressure
over some coarse-grid blocks.

\section{Multiscale inversion} \label{sec:method}

In this section, we describe the inversion formulation, and the numerical algorithm.
\subsection{Inversion formulation}
Denote the fine-grid solution by $u_h$ and the coarse-grid solution by
\[
u_H=\sum_{i,j} c_{ij}\phi_j^{\omega_i},
\]
where $\phi_j^{\omega_i}(x)$ are GMsFEM basis functions, which can approximate the fine-grid solution $u_h$ accurately
for the inverse problem. We shall denote the vector of expansion coefficients $c_{ij}$ by $c$. Throughout, it is always assumed
that the problem has a reduced dimensional approximation, i.e., very few basis functions can provide a good approximation
of the fine-grid solution $u_h$ (in a suitable norm $\|\cdot\|$):
\begin{equation}\label{eqn:prox}
\|u_h-u_H\|\approx \text{small}.
\end{equation}

Suppose that we measure the quantity $F_{\text{obs}}$ defined by
\[
F_{\text{obs}}=G(u_h),
\]
where $G$ is a bounded linear function. In view of the relation \eqref{eqn:prox}, we have $G(u_h)\approx G(u_H).$
Next, we formulate the inverse problem in terms of discrete parameters (defined on the coarse grid). Note that the
coefficient vector $c$ of the discrete coarse-grid solution $u_H$ has a form
\[
M{\frac{dc}{dt}} + A c= b,
\]
with unknown low dimensional matrices $A$ and $M$ (which depend on basis functions $\phi_j^{\omega_i}(x)$ and $\kappa$
-- both are unknown in the inverse context), and the time-dependent vector $b$ is source term. By the linearity of the
operator $G$, we also have
\[
F_{\text{obs}}\approx G(u_H)=\sum_{i,j} c_{i,j} G(\phi_j^{\omega_i}).
\]

For the proposed multiscale inversion technique, the standing assumption on the measurement operator $G$ is that
\begin{equation}\label{assump1}
G(\phi_j^{\omega_i})=y(c,A,M),
\end{equation}
i.e.,  the observed response $F_{\text{obs}}$ can be expressed in terms of the elements of the stiffness and mass matrices $A$ and $M$
and the coefficient vector $c$. This assumption holds true for a wide variety of observations, which are averaged quantities
over coarse blocks, e.g., pressures or fluxes. In this case, we have
\[
F_{\text{obs}}=\mathcal{Y}(c,A,M).
\]

We illustrate this general formulation with two more concrete examples.
For example, if we observe the average pressure on a coarse block $K$ away from the boundary:
\[
{y}_K=\int_K  u_H {\rm d}x= c_{ij} \int_K \phi_i^{\omega_j}{\rm d}x.
\]
To express the given data this in terms of $c$, $A$, and $M$, we recall the entry $(M)_{ij,kl}$ of the mass matrix
$(M)_{ij,kl}=\int_\omega \phi_i^{\omega_j} \phi_k^{\omega_l}{\rm d}x$. In our numerical studies, we
seek the element-wise components of $(M)_{ij,kl}$ for each $K$ (see Fig. \ref{fig:illustration}),
denote it by $(M)^K_{ij,kl}$. Then, since the first basis functions form the partition of unity, there holds
\[
{y}_K = c_{ij}\sum_{k=1,l\in \mathcal{I}} (M)_{ij,kl}^K,
\qquad G=\mathcal{G}(c,M),
\]
where $\mathcal{I}$ is the set of indices for coarse vertices. Similarly,  if we observe the average flux
(for simplicity, we denote it by $y_K$) over a coarse block $K$
\[
{y}_K=\int_K \kappa \nabla u_H {\rm d}x= c_{ij} \int_K \kappa \nabla \phi_i^{\omega_j}{\rm d}x.
\]
Note that
$
(A)_{ij,kl}=\int_\omega \kappa \nabla\phi_i^{\omega_j} \cdot \nabla \phi_k^{\omega_l} {\rm d}x.
$
Then, one can solve for $\int_K \kappa \nabla u_H{\rm d}x$ from
$(A)_{ij,kl}^K$, the elements of the stiffness matrix in $K$
corresponding to $k=1$. To do so, we first note that
$
(A)_{ij,1l}^K=\int_K \kappa \nabla\phi_i^{\omega_j} \cdot \nabla \phi_1^{\omega_l} {\rm d}x=\int_K \kappa \nabla\phi_i^{\omega_j} \cdot \nabla \phi_{\omega_l}^0 {\rm d}x,
$
where $ \phi_{\omega_l}^0$ are linear basis functions.
By solving the resulting $2\times 2$ system, we can compute
$\int_K \kappa \nabla\phi_i^{\omega_j}{\rm d}x$.


Note that in this case, we cannot identify the solution $u(x)$ explicitly, since we do not
know the basis functions. However, given the elements of the stiffness matrix $A$, we can
find some properties of the fine-grid permeability $\kappa(x)$. Upon writing the observation
in terms of coarse-grid discretization parameters, the multiscale inverse problem has
the following formulation
\begin{equation}
P(c,A,M|d)\propto P(d|c,A,M) \pi(A) \pi(M) \pi(c).
\end{equation}
The priors on $A$, $M$, and $c$ can be specified in various ways. In our simulations, we use Gaussian
priors around a given state generated with a fixed permeability field. In general, one can use a Gaussian mixture
field based on several generated permeability fields or priors generated using fine-grid permeability
fields as in a Bayesian framework \cite{ChenZabaras:2015}; however, we stress that our objective is to recover coarse-grid parameters.
Once we identify $c$, $A$, and $M$, some solution averages can be obtained. To formalize this process,
we assume that we can construct a set of observable basis functions $\widetilde{\phi_j^{\omega_i}}$ such that
\[
\sum_{i,j} c_{i,j} G(\phi_j^{\omega_i})\approx \sum_{i,j} c_{i,j} G(\widetilde{\phi_j^{\omega_i}}),
\]
or equivalently $G(\phi_j^{\omega_i})\approx G(\widetilde{\phi_j^{\omega_i}}).$
This latter assumption has to be verified case by case for each operator $G$. Generally, it is necessary for performing
inversion on a coarse grid in order to guarantee that the observation can be observed on a coarse-grid solution.

\begin{remark}
When the permeability field $\kappa(x)$ is parameterized or samples of permeability fields are known, we can compute
the multiscale basis functions $\phi_i^{\omega_j}$. Then one can compute the fine-grid solution without explicitly finding $\kappa$.
\end{remark}

\begin{remark}[One basis function - numerical homogenization]
In numerical homogenization, our goal is to find $\kappa^*$ on a coarse grid. Then the coarse-grid solution $u_H$ satisfies
\begin{equation*}
{\frac{\partial u_H^*}{\partial t}} + \mathcal{L}(\kappa^*,u_H^*)=0,
\end{equation*}
where $\mathcal{L}$ is an elliptic differential operator depending on $\kappa^*$. Assume that we can observe the data $F_{obs}$ based on a
coarse-grid solution $u_H$: $G(u^*_H)=F_{\text{obs}}$. In analogy, we assume that one un-observable basis
function can be used to approximate the solution $u_H= \sum_i c_i \phi^{\omega_i}.$ Then we can take
$\widetilde{\phi^{\omega_i}}=\phi_0^{\omega_i}$, polynomial basis function that has the same linear boundary
conditions as multiscale basis functions.
\end{remark}

\begin{remark}[Multi-continuum approach]
In the recent work \cite{ChungEfendievLeungVasilyeva:2017}, we have discussed the relation between the GMsFEM
and multi-continuum approaches. For multi-continuum equations, the generalized multiscale inversion
technique reduces to finding parameters in multi-continuum equations. For example, in a simplified case,
the coarse-grid equations assume the form
\[
{\frac{\partial u_{i,H}^*}{\partial t}} - div(\kappa_i^*\nabla u_{i,H}^*) +
Q_{ij}(u_{j,H}^* -u_{i,H}^*)=0,
\]
where the index $i$ refers to the continua and our goal is to identify
$\kappa_i^*$ and $Q_{ij}$. The latter can be done using
standard inverse problem approaches. As a result, we compute
the effective properties of each continua and they cannot be directly
related to the fine-grid permeability field. Our approach
avoids assumptions of multi-continua approaches and, while as
in multi-continua inversion, identifies some average properties
about the fine-scale permeability field.

\end{remark}

\subsection{Numerical algorithm}

In practice, the inversion can be performed by solving the following minimization problem
\begin{equation}\label{eq:sigmas}
J(M, A, u) =\frac{1}{\sigma_M^2} ||M - M_0||^2 + \frac{1}{\sigma_A^2} ||A - A_0||^2 + \frac{1}{\sigma_F^2} ||Fu - g||_{L^2(0, T)}^2,
\end{equation}
where $M$ and $A$ are global mass and stiffness matrices, respectively, and $M_0$ and $A_0$ are the corresponding
given prior information. Below, we use Gaussian priors around a state generated with a fixed permeability field.
In general, one can use a Gaussian mixture model based on several generated permeability fields or priors generated by
fine-grid permeability fields as in Bayesian models, as mentioned earlier. We remark that the positive scalars $\sigma_M, \sigma_A$ and
$\sigma_F$ play the role of regularization parameters, which are constructed to give relative weights of the terms.
Choosing proper regularization parameters is a notoriously challenging task in general and depends on the choice of
the prior and the specific application, and we refer to  \cite{ItoJin:2015} for detailed discussions on various ways for selecting a single
regularization parameter. Moreover, one needs to select the norms appropriately in \eqref{eq:sigmas} to guarantee
the robustness with respect to the data perturbation \cite{ItoJin:2015}. In this work, for simplicity, we employ the discrete
quantities and $l^2$ norms in \eqref{eq:sigmas}, and leave the rigorous studies to a future work. In the functional,
the operator $F$ is the measurement operator, and $g$ is the observed data.
In our numerical simulations, the observed data $g$ is obtained by solving the forward problem on the fine grid,
and then apply the operator $F$ to the solution $u_h$. In particular, for each coarse element $K$, we have
\begin{equation*}
F^K u_h:= g^K(t) := \bar{u}_h^K(t) = \frac{1}{|K|} \sum_{i = 1}^{DOF_K} c^K_i(t) \sum_{j=1}^{DOF_K} m^K_{ij}.
\end{equation*}

We will solve the optimization problem \eqref{eq:sigmas} using an iterative procedure. First, we assume that some
initial approximations for the local mass and stiffness matrices $A_0^K$ and $M_0^K$ are given. These matrices are
found by generating a priori realization and used as a regularization for the low-dimensional inverse problem. Based
on these initial conditions, we solve the global forward problem and find an initial approximation $u_0(t)$
\begin{equation*}
  (A_0, \quad M_0) \rightarrow  u_0(t),
\end{equation*}
and the corresponding simulated observational data $g_0^K(t)$ is the average pressure in cell $K$
\begin{equation*}
g_0^K(t) = \bar{u}_0^K(t)
= \frac{1}{|K|} \sum_{i = 1}^{DOF_K} c^K_{0,i} (t) \sum_{j=1}^{DOF_K} m^K_{ij},
\end{equation*}

The multiscale inversion algorithm proceeds as follows. First, we seek the element-wise components of the
stiffness and mass matrices $A$ and $M$. In this way, we can ensure the symmetry. We iteratively ($n = 1,
2, ...$) update local mass matrix $M_n^K$ and local stiffness matrix $A_n^K$
\begin{equation}\label{eq:sys1}
M^K_n = M^K_{n-1}  - \epsilon \delta J_M \quad \mbox{and}\quad A^K_n = A^K_{n-1}  - \epsilon  \delta J_A,
\end{equation}
using the previous iterates $M^K_{n-1}$ and $A^K_{n-1}$, where $\epsilon >0$ is the step size, and $\delta J_A$
and $\delta J_M$ denote the derivative of the functional $J$ with respect to $A$ and $M$, respectively.
Further, we generate global mass and stiffness matrices by local matrices and solve the global forward problem
\begin{equation}
(A_n, \quad M_n) \rightarrow  u_n(t),
\end{equation}
and find average cell pressure
\begin{equation}
g_n^K(t) = \frac{1}{|K|} \sum_{i = 1}^{DOF_K} c^K_{n,i}(t) \sum_{j=1}^{DOF_K} m^K_{ij}.
\end{equation}

In the gradient descent iteration \eqref{eq:sys1}, we need the derivatives $\delta J_M$ and $\delta J_A$ of the
functional $J$ with respect to the matrices $M$ and $A$. To this end, we employ the standard adjoint state technique.
In the following, we only give the main steps since the derivation is rather standard \cite{ItoKunisch:2008}.
Consider the adjoint problem
\[
{\frac{\partial w}{\partial t}} + \text{div} (\kappa\nabla w) =-F^T(Fu-g), \quad w(T)=0
\]
where $F^T$ is the adjoint of the operator $F$. Note that the adjoint problem is defined backward in time,
and can be solved numerically as usual by a change of variable $t\leftarrow T-t$. Suppose that we represent
the adjoint solution $w$ as $\{ \lambda_{n-1,i} \}$ in the multiscale basis $\phi_i^{\omega_j}$. Then using
the adjoint solution $w(t)$ and the forward solution $u_{n-1}(t)$, the local component
$(\delta J_M)^K_{ij}$ of the derivative $\delta J_M$ can be computed as
\begin{equation*}
(\delta J_M)^K_{ij} = \frac{2}{\sigma_M^2} \Big( (M_{n-1})^K_{ij} - (M_0)^K_{ij} \Big)
- \frac{2}{\sigma_F^2} (M_{n-1})^K_{ij} \int_0^T  \frac{d\lambda_{n-1,j_g}}{dt} c_{n-1,i_g},
\end{equation*}
where $i_g$ is the corresponding global index. That is, $i_g$ is the global index of the vertex corresponding to the
local index $i$. Similarly, we can compute the derivative $\delta J_A$ as
\begin{equation*}
(\delta J_A)^K_{ij} = \frac{2}{\sigma_A^2} \Big( (A_{n-1})^K_{ij} - (A_0)^K_{ij} \Big)
- \frac{2}{\sigma_F^2} (A_{n-1})^K_{ij} \int_0^T  \lambda_{n-1,j_g} c_{n-1,i_g}.
\end{equation*}


\section{Numerical results}\label{sec:num}
Now we illustrate our multiscale inversion technique with flows in fractured media, where fractures have
high conductivities and very small width. In our simulations, their widths are assumed to be zero and they are
modeled as high-conductivity lines; see \cite{ChungEfendievLeungVasilyeva:2017} for details. The presence of
multiple disconnected fracture networks requires several basis functions as we discussed earlier.

\begin{figure}[!ht]
\centering
\includegraphics[width=0.99 \textwidth]{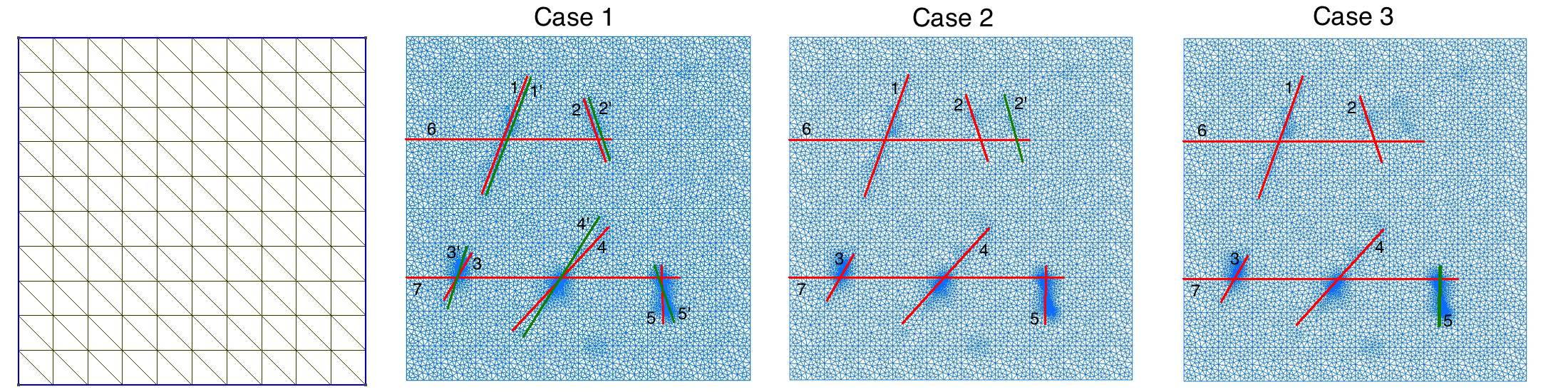}
\caption{Coarse and fine grids for Cases 1-3. In the figure, red color indicates exact fractures, and the
green color is for moved fractures. Case 1 has 3 rotated and 2 shifted  fractures,
Case 2 has 1 shifted fracture with large distance (second) and  Case 3 has 1 removed fracture (fifth).
}
\label{fig:mesh}
\end{figure}

\begin{figure}[!ht]
\centering\setlength{\tabcolsep}{0mm}
\begin{tabular}{ccc}
\includegraphics[width=0.32 \textwidth]{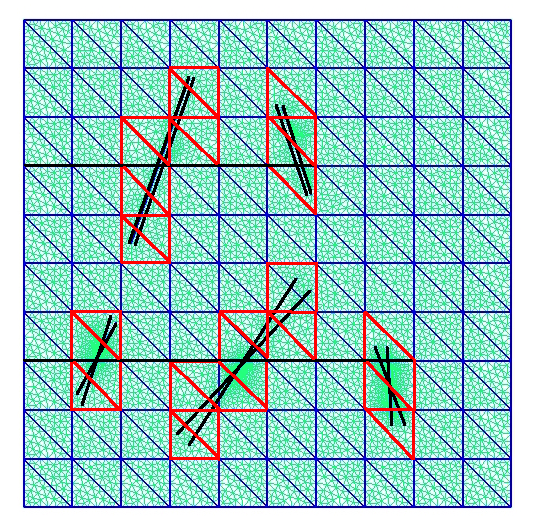}&
\includegraphics[width=0.32 \textwidth]{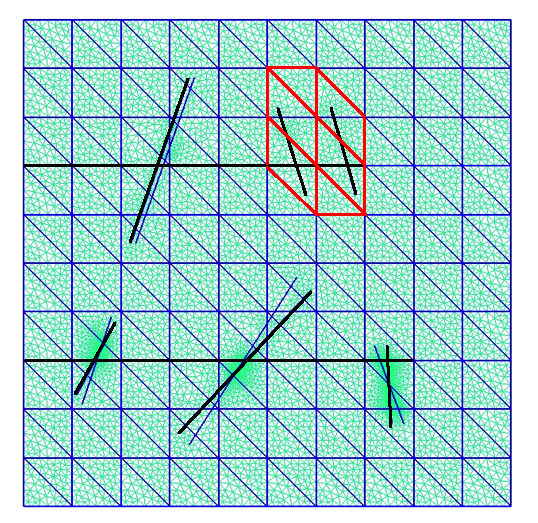}&
\includegraphics[width=0.32 \textwidth]{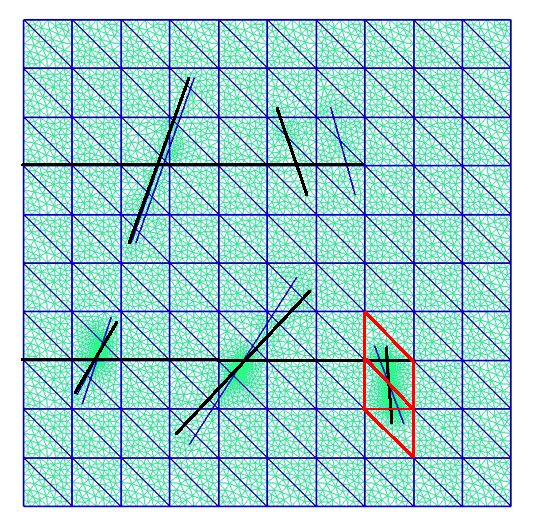}\\
(a) case 1 & (b) case 2 & (c) case 3
\end{tabular}
\caption{The regions for (adaptive) inversion update for the three cases: The coarse cells in red indicate the corresponding
entries of the matrices to be updated.} 
\label{fig:mesh-adapt}
\end{figure}

In our numerical experiment, we take the computational domain $\Omega = [0,1]^2$. The coarse mesh contains 121
vertices and 200 cells. We use the following parameters
\begin{itemize}
\item $k_m = 10^{-3}$ and $k_f = 10^{2}$,
\item $c_m = c_f = 1.0$,
\item $p = 0$ on the left boundary and no flow on the remaining boundaries with $p_0 = 1$ for $t = 0$,
\item $f = 0$ and $t_{max} = 10$ with 10 time steps.
\end{itemize}
The fine mesh contains 6297 vertices and 12352 cells for Case 1. For Case 2, we have 7908 vertices and 15574 cells.
For Case 3, fine mesh with 7891 vertices and 15540 cells. The fine meshes for all three cases are depicted in Fig. \ref{fig:mesh}.
In Fig. \ref{fig:mesh-adapt}, we show the adaptive regions, where we perform updates to the matrices, and unless
otherwise stated, these regions are used in all the numerical experiments with the proposed inversion technique below.
Further, unless otherwise stated, the step length $\epsilon$ in the iteration \eqref{eq:sys1}
is fixed at $\epsilon=10^{-12}$.

\begin{figure}[hbt!]
\centering
\begin{tabular}{cc}
\includegraphics[width=0.5 \textwidth]{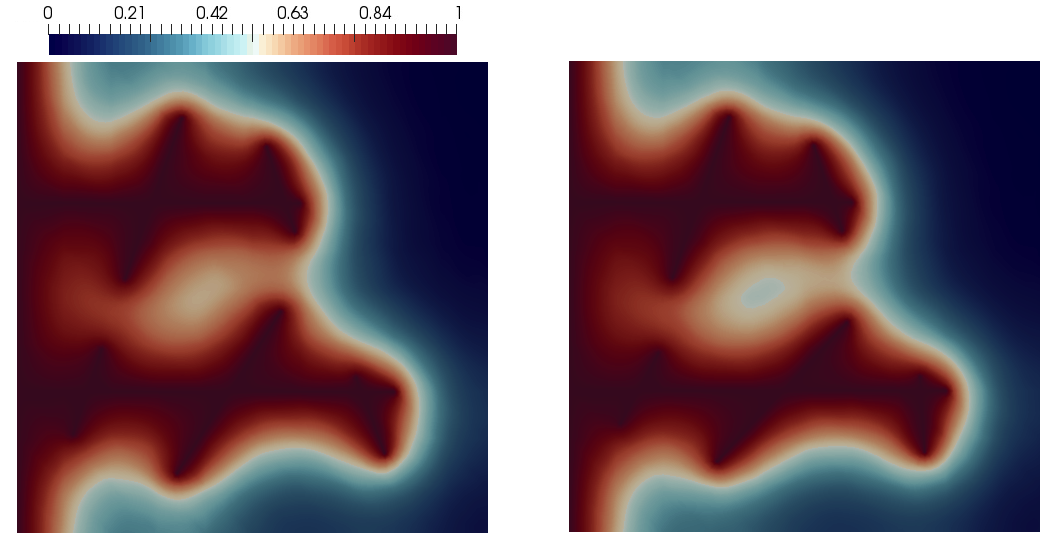}&
\includegraphics[width=0.5 \textwidth]{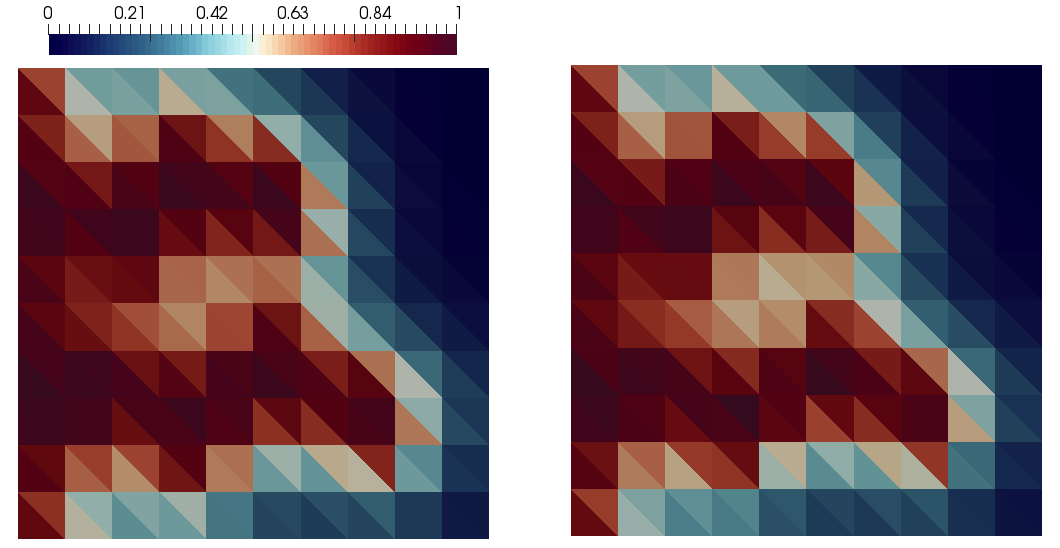}\\
(a) multiscale solution on fine grid & (b) 4 multiscale basis functions, adaptive inverse
\end{tabular}
\caption{Numerical results Case 1:
(a) multiscale solution for $u_0$ (left)  and exact (right), and (b) cell average solution for initial condition $M_0$, $A_0$
(left) and solution after 100 iterations (right).}
\label{fig:inv3-u}
\end{figure}

\begin{figure}[hbt!]
\centering
\begin{tabular}{cc}
\includegraphics[width=0.5 \textwidth]{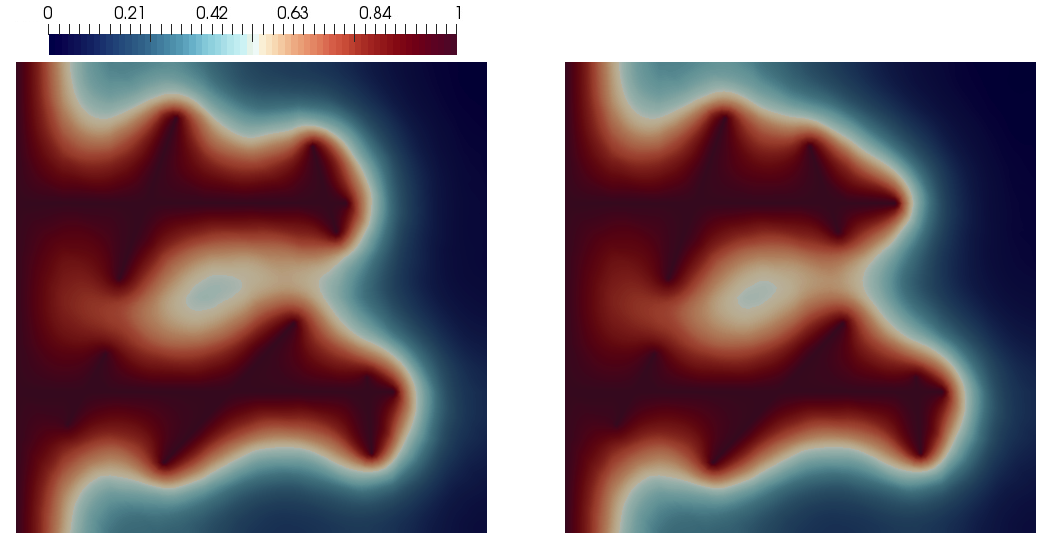} & \includegraphics[width=0.5 \textwidth]{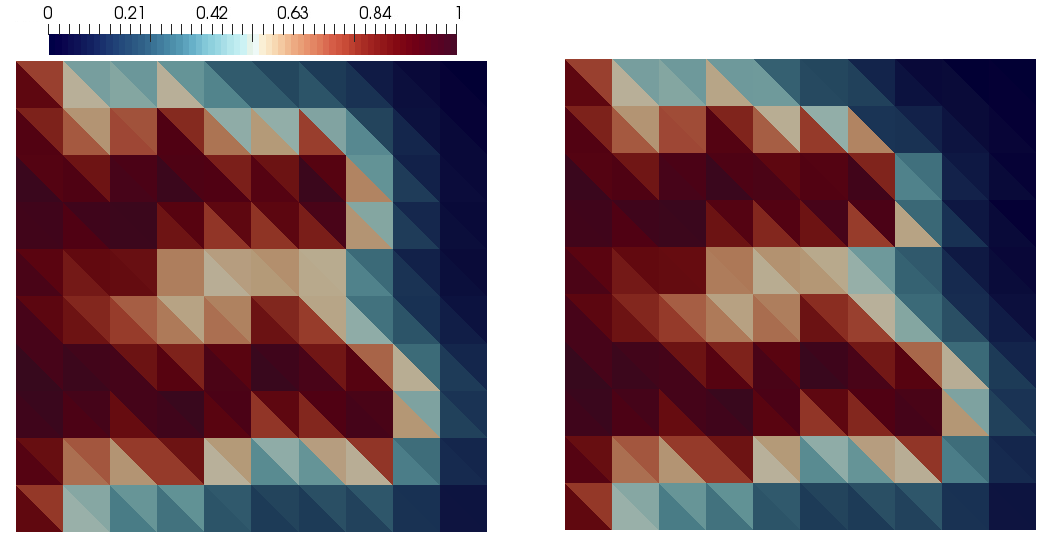}\\
(a) multiscale solution on fine grid         & (b) 4 multiscale basis functions, adaptive inverse\\
\end{tabular}
\caption{Numerical results for Case 2: (a) multiscale solution for $u_0$ (left)  and exact (right), and
(b) cell average solution for initial condition $M_0$, $A_0$  (left)  and solution after 100 iterations (right).}
\label{fig:inv4r-u}
\end{figure}

\begin{figure}[hbt!]
\centering
\begin{tabular}{cc}
\includegraphics[width=0.5 \textwidth]{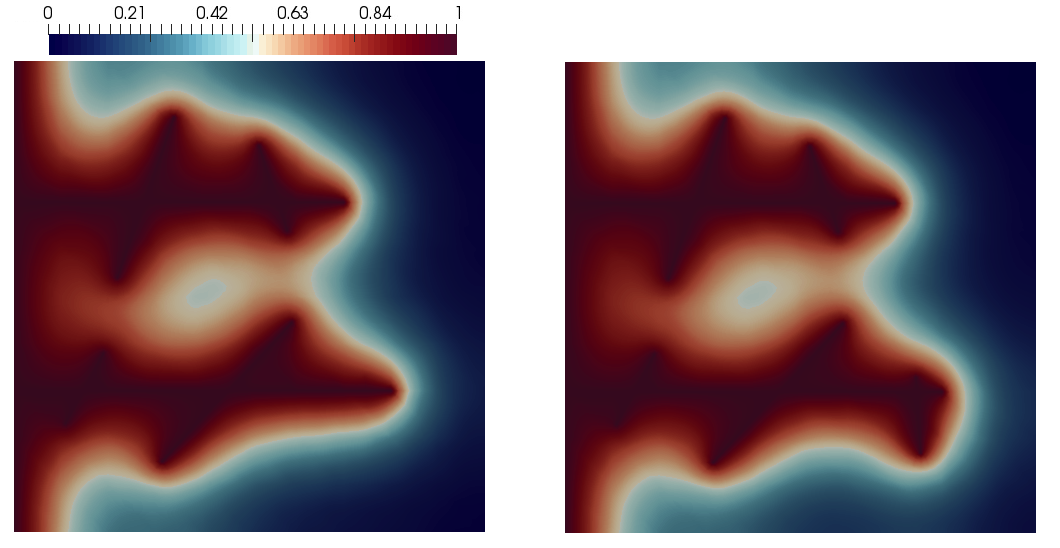} & \includegraphics[width=0.5 \textwidth]{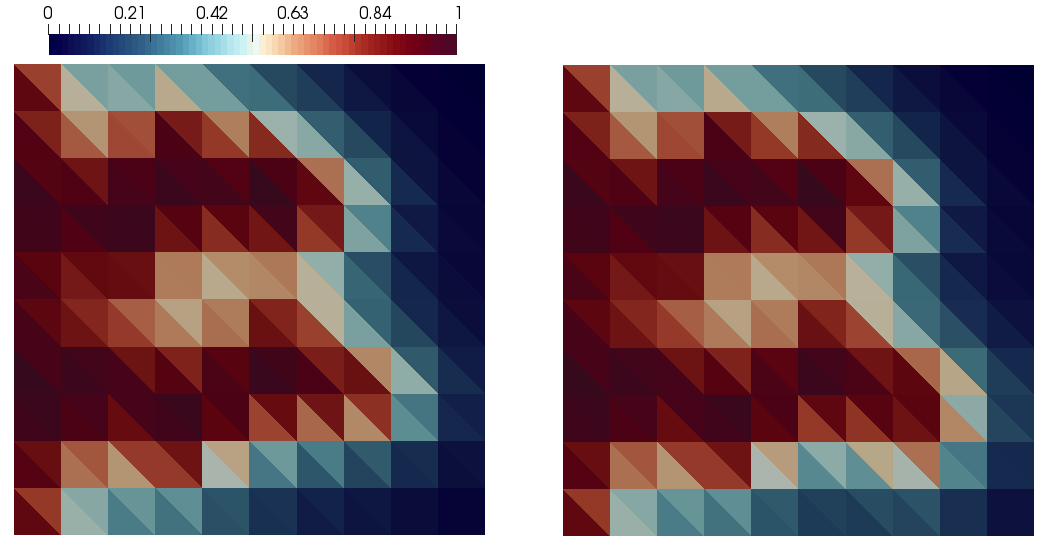}\\
(a) multiscale solution on fine grid & (b) 4 multiscale basis functions, adaptive inverse
\end{tabular}
\caption{Numerical results for Case 3: (a) multiscale solution for $u_0$ (left)  and exact (right), and (b)
cell average solution for initial condition $M_0$, $A_0$  (left) and solution after 100 iterations (right).}
\label{fig:inv5r-u}
\end{figure}

To evaluate the proposed approach, we first present the following numerical results: the average fine-grid solution,
the initial condition and the final solution. All the results are obtained with the following parameter setting:
$\sigma_M = \sigma_A = 1.0$ and $\sigma_F = 10^4$, which are determined in a trial and error manner. In Figs.
\ref{fig:inv3-u}, \ref{fig:inv4r-u} and \ref{fig:inv5r-u} for the three cases, we present the numerical solutions. It is observed
that the recovered solutions are always fairly close to the exact one, indicating the accuracy of the proposed approach.

\begin{figure}[hbt!]
\centering
\begin{tabular}{cc}
\includegraphics[width=0.4 \textwidth]{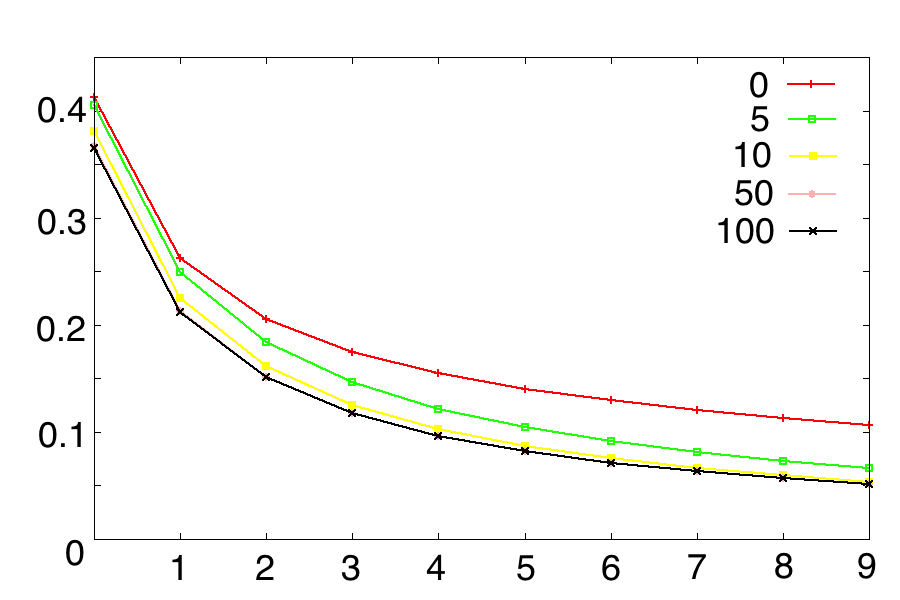}
& \includegraphics[width=0.4 \textwidth]{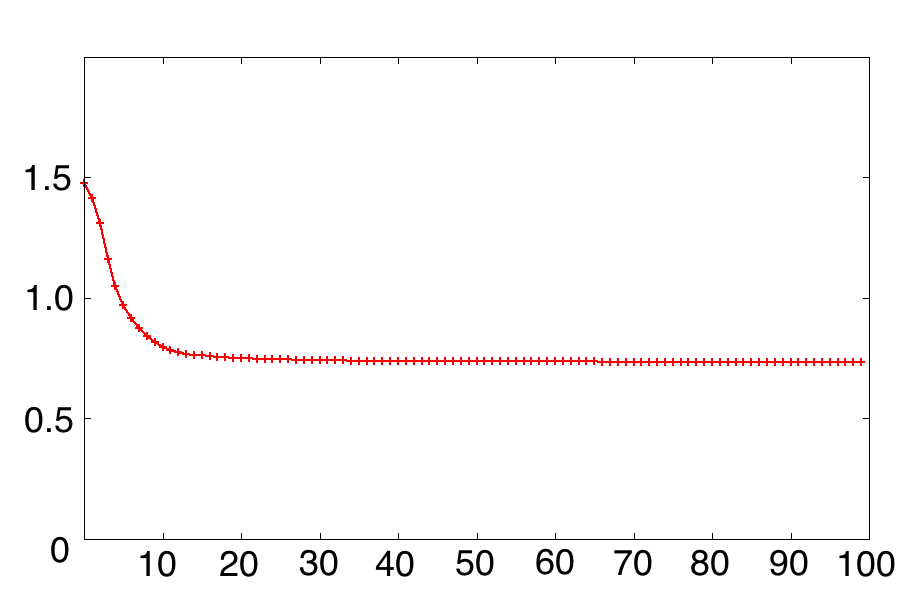}\\
 \includegraphics[width=0.4 \textwidth]{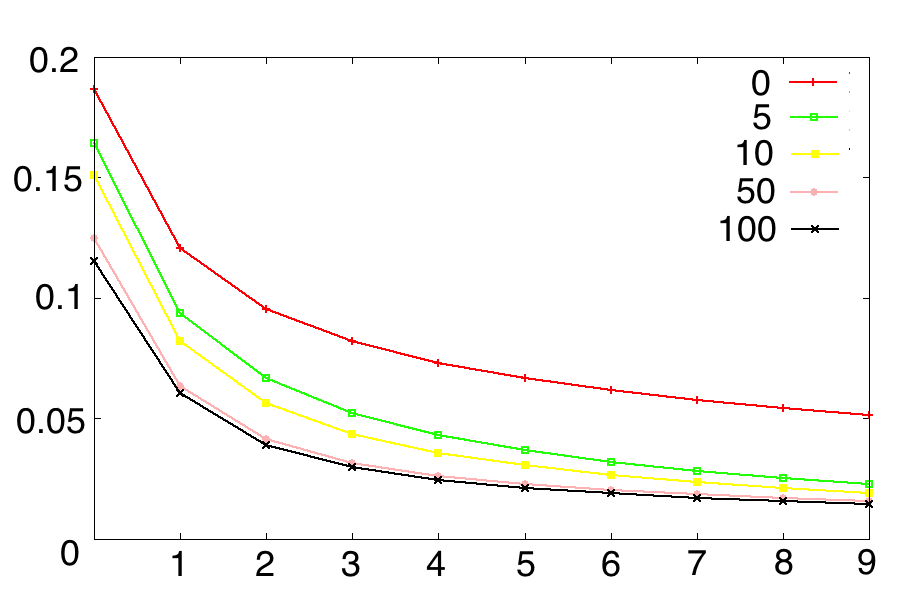}&
\includegraphics[width=0.4 \textwidth]{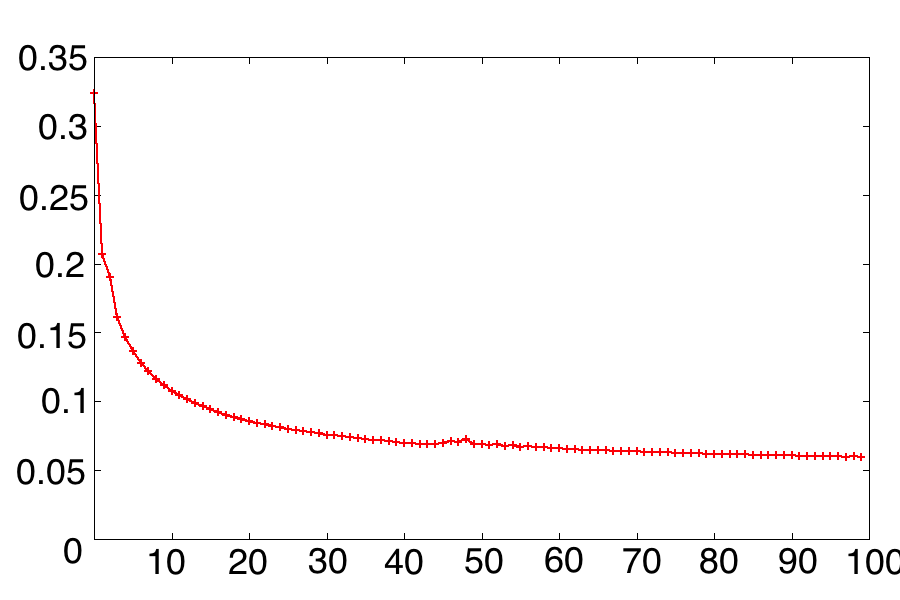}\\
(a) $L^2$ error v.s. time & (b) $J$ v.s. iteration
\end{tabular}
\caption{Numerical results for Case 1 with 2 (top) or 4 (bottom) multiscale basis functions:
(a) the $L^2$ error for cell average v.s. time, and (b) the functional value
$J$ v.s. iteration index.}\label{fig:inv3-e1}
\end{figure}
\begin{figure}[hbt!]
\centering
\begin{tabular}{cc}
\includegraphics[width=0.4 \textwidth]{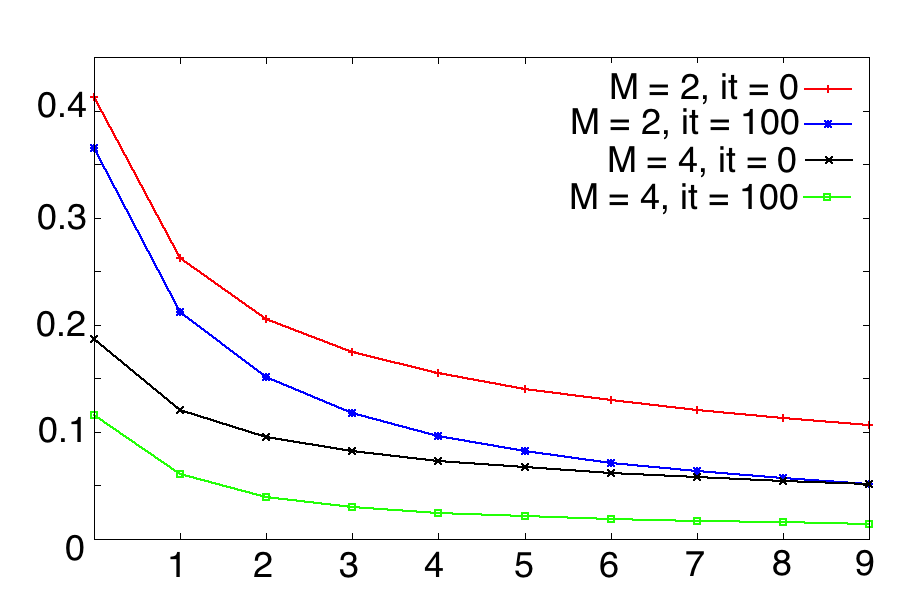}&\includegraphics[width=0.4 \textwidth]{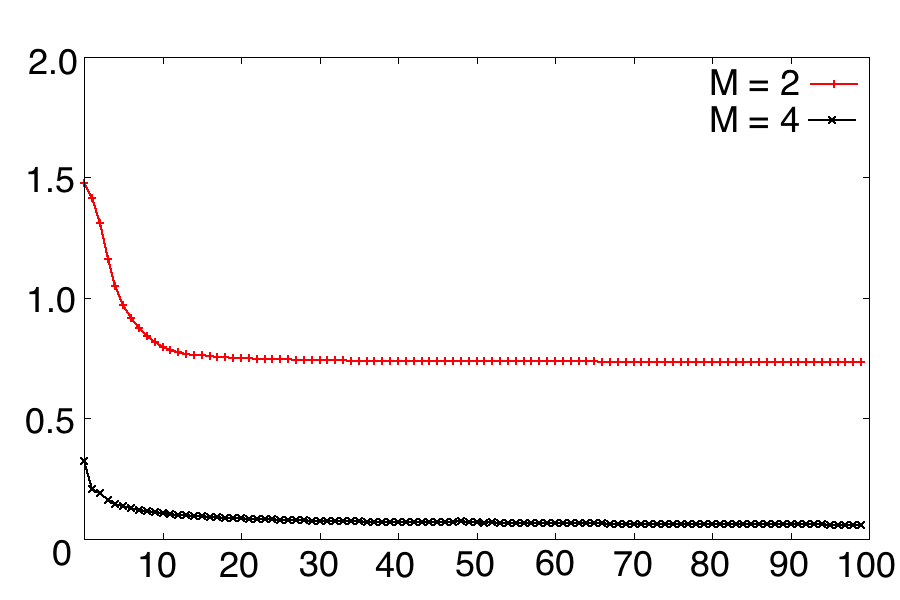}\\
(a) $L^2$ error v.s. time & (b) $J$ v.s. iteration
\end{tabular}
\caption{Numerical results for Case 1, using 2 or 4 multiscale basis functions: (a) $L^2$ error for cell average v.s. time $t$, and (b)
the functional value $J$ v.s. the iteration index.}
\label{fig:inv3-e}
\end{figure}
\begin{figure}[hbt!]
\centering
\begin{tabular}{cc}
\includegraphics[width=0.4 \textwidth]{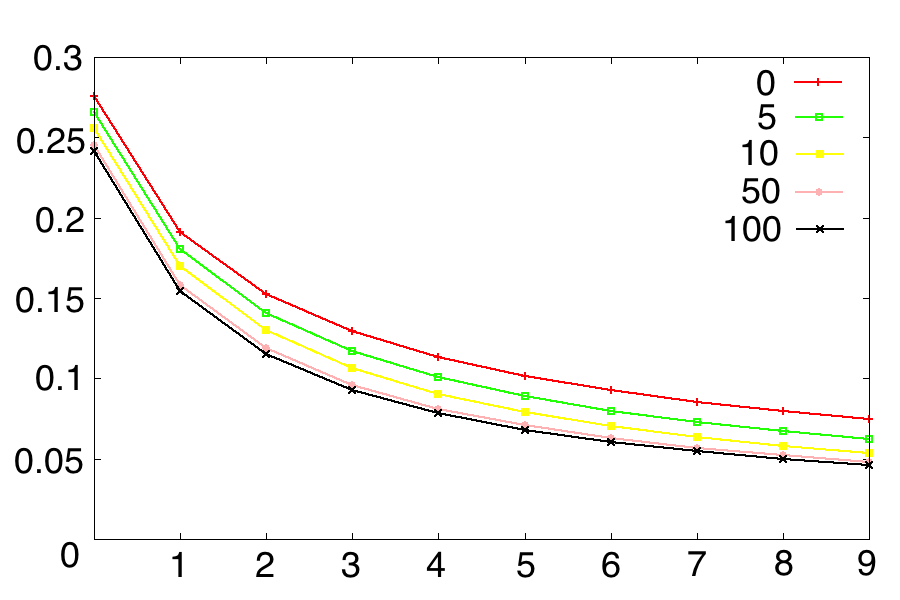}&
\includegraphics[width=0.4 \textwidth]{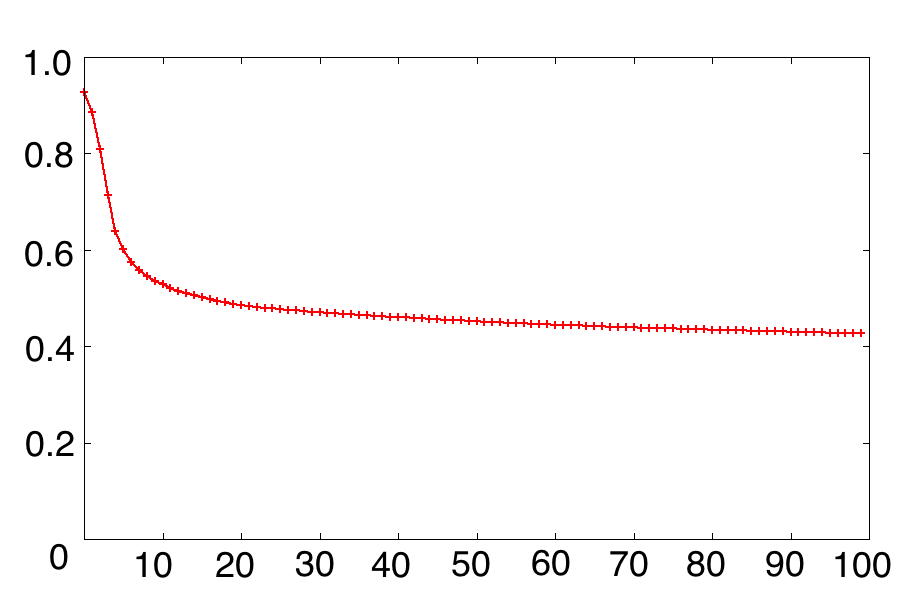}\\
(a) $L^2$ error v.s. time & (b) $J$ v.s. iteration
\end{tabular}
\caption{Numerical results for Case 2, using 4 multiscale basis functions: (a) $L^2$ error for cell average v.s. time $t$ at different iterations, and (b)
the functional value $J$ v.s. the iteration index.}
\label{fig:inv4r-4b}
\end{figure}
\begin{figure}[hbt!]
\centering
\begin{tabular}{cc}
\includegraphics[width=0.4 \textwidth]{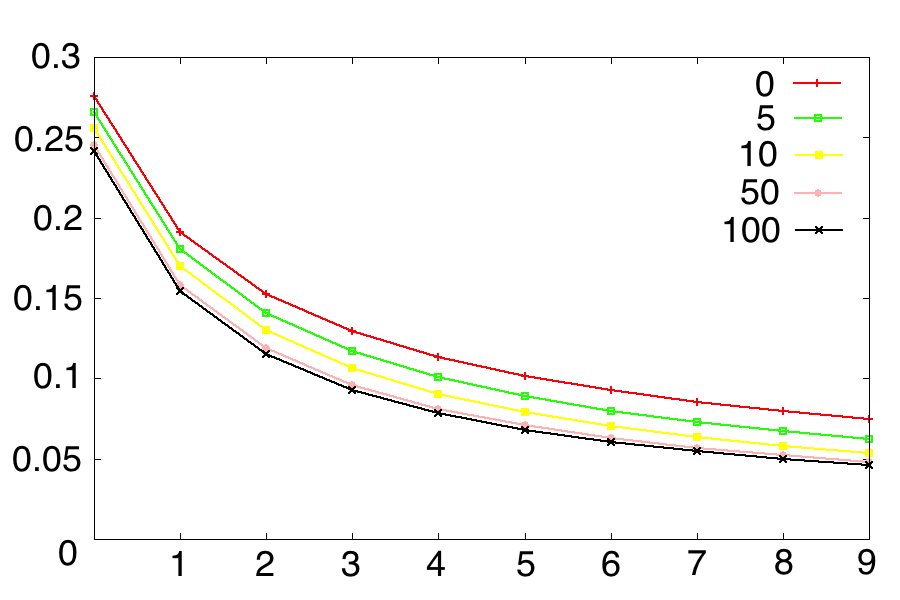}&
\includegraphics[width=0.4 \textwidth]{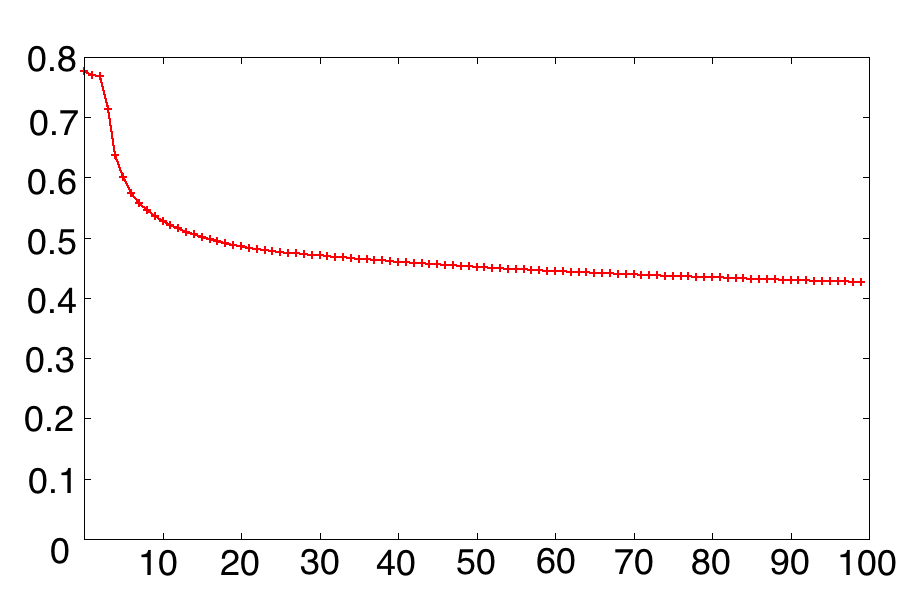}\\
(a) $L^2$ error v.s. time & (b) $J$ v.s. iteration
\end{tabular}
\caption{Numerical results for Case 3 using 4 multiscale basis functions: (a)
$L^2$ error for cell average v.s. time $t$ for different iterations, and (b) the functional value $J$ v.s. the iteration index.}
\label{fig:inv5r-4b}
\end{figure}

In Figs. \ref{fig:inv3-e1}, \ref{fig:inv3-e}, \ref{fig:inv4r-4b} and \ref{fig:inv5r-4b}, we present the $L^2$
errors with respect to the space variable as a function of time $t$ and the residual (functional value) $J$ given in \eqref{eq:sigmas}.
The $L^2$ error decreases with the time $t$, and in the absence of data noise, it also decreases as the
iteration proceeds. Further, with more multiscale basis functions in the inversion, the $L^2$ error is also smaller.
We always observe that the residual decreases as the number of iterations grows. The monotone decrease of the residual
indicates that the iteration \eqref{eq:sys1} is indeed minimizing the functional $J$.

\begin{figure}[hbt!]
\centering
\begin{tabular}{cc}
\includegraphics[width=0.4 \textwidth]{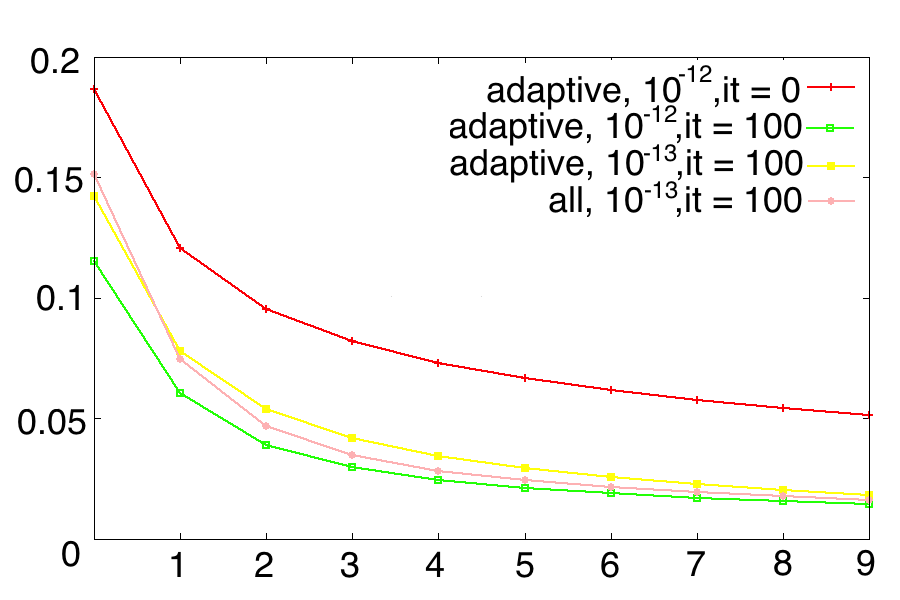}&
\includegraphics[width=0.4 \textwidth]{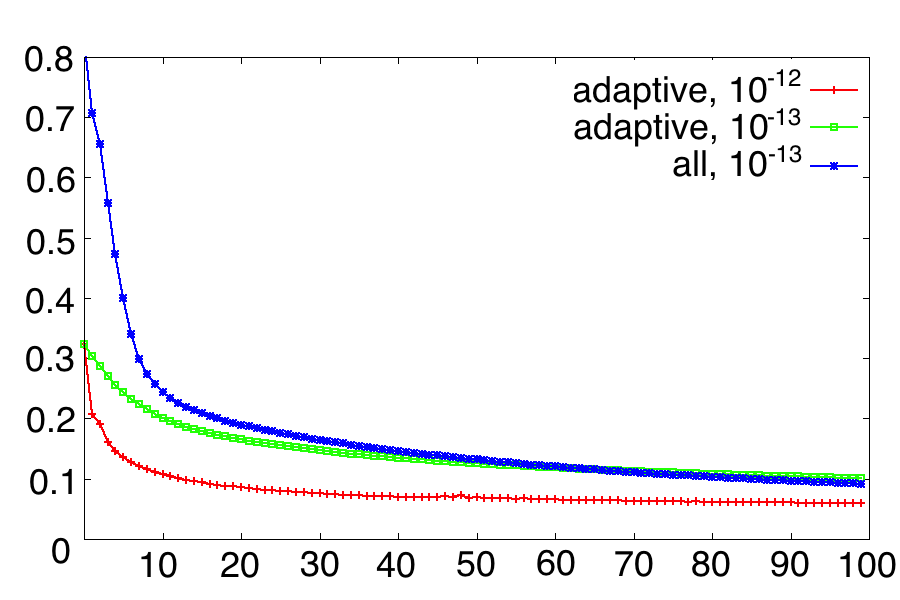}\\
(a) $L^2$ error v.s. time & (b) $J$ v.s. iteration
\end{tabular}
\caption{Numerical results for Case 1, using 4 multiscale basis functions:
(a) $L^2$ error for cell average v.s. time $t$, and (b) the functional value $J$ v.s. iteration index.
Different iteration parameter $\epsilon = 10^{-12}$ and $10^{-13}$. Adaptive and all cells local mass matrices updating. }
\label{fig:inv3-ad}
\end{figure}

\begin{figure}[hbt!]
\centering
\begin{tabular}{cc}
\includegraphics[width=0.4 \textwidth]{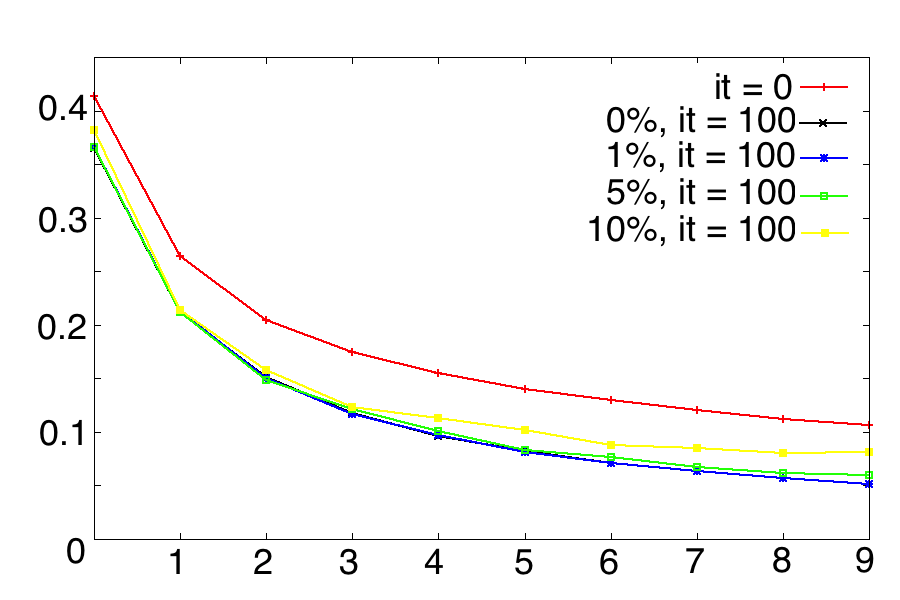}&
\includegraphics[width=0.4 \textwidth]{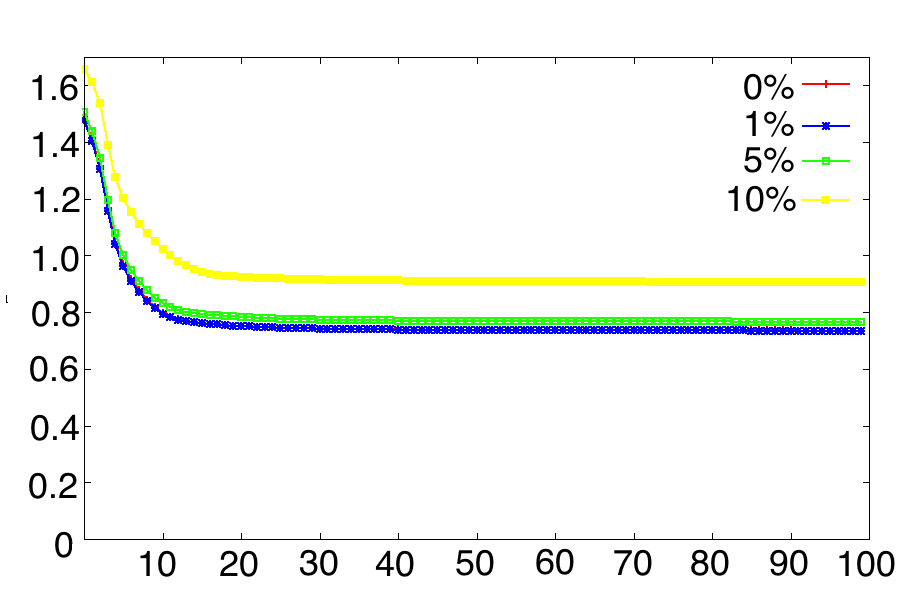}\\
\includegraphics[width=0.4 \textwidth]{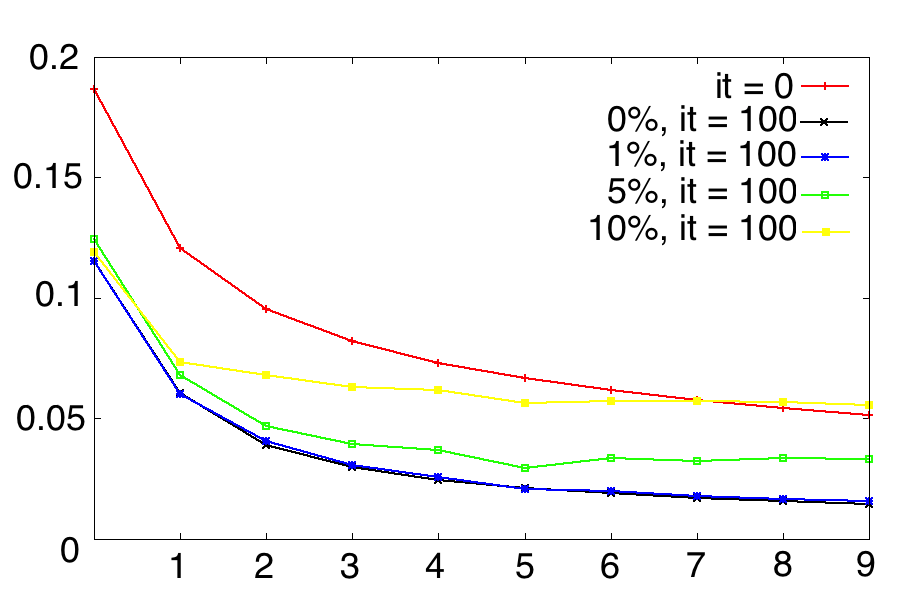}&
\includegraphics[width=0.4 \textwidth]{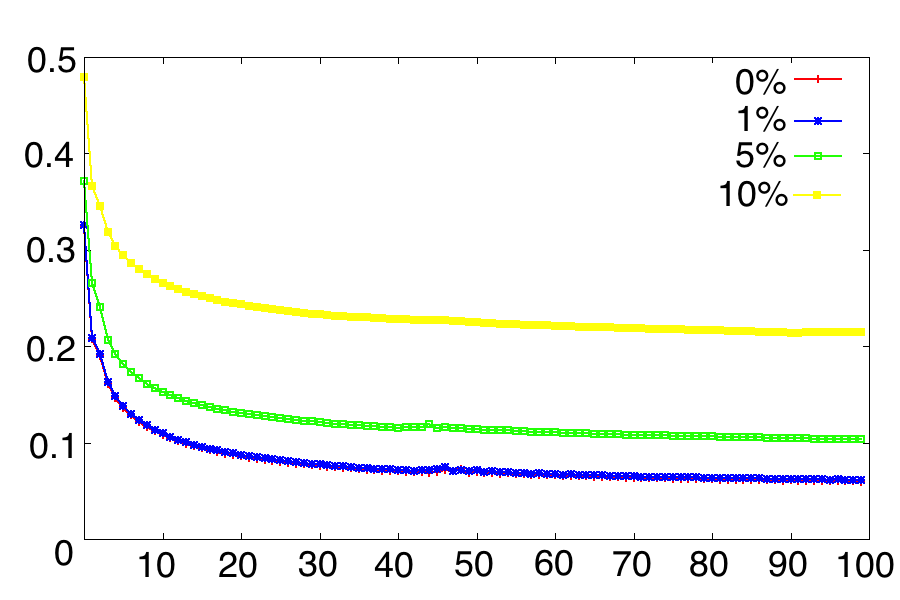}\\
(a) $L^2$ error v.s. time & (b) $J$ v.s. iteration
\end{tabular}
\caption{Numerical results for Case 1 with noisy data $g$, $g^K(t) = (1 + \delta r) g^K(t) $ ($r \in [-1, 1]$ is random number
and $\delta = 1 \%$, $3 \%$ , $5 \%$ or $10 \%$) with 2 (top) or 4 (bottom) multiscale basis functions:
(a) the $L^2$ error for cell average v.s. time, and (b) the functional value $J$ v.s. iteration index.}
\label{fig:inv3r-rand}
\end{figure}

\begin{figure}[hbt!]
\centering
\includegraphics[width=0.9 \textwidth]{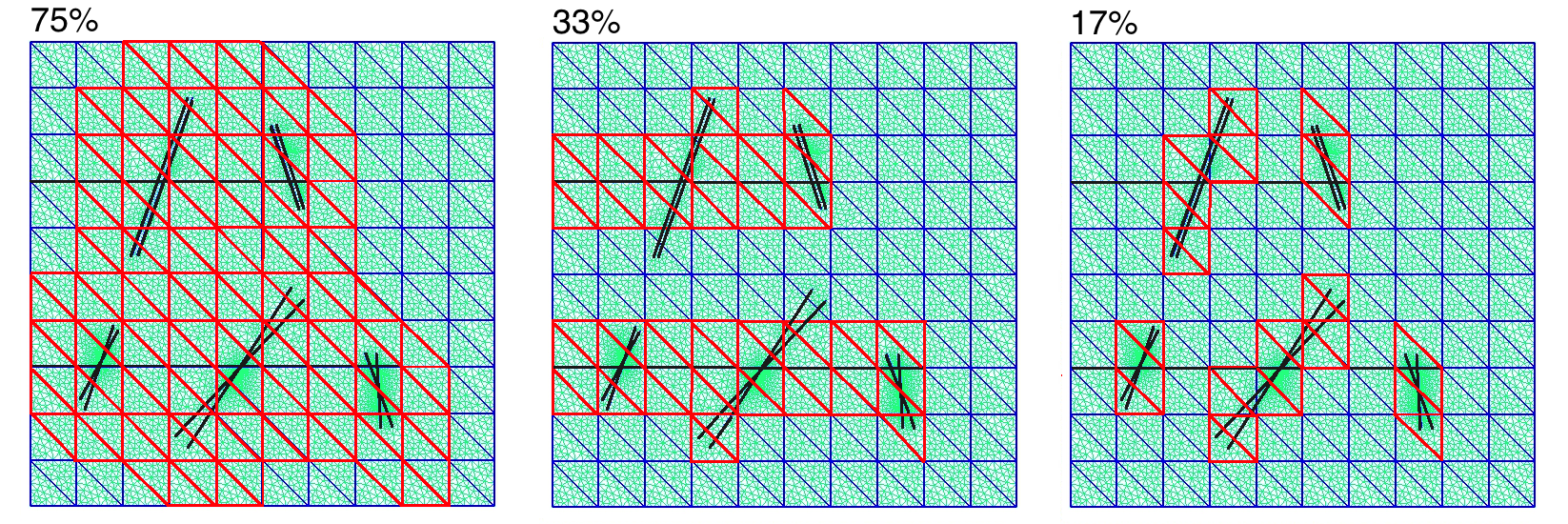}
\caption{The observation data $g^K$ for Case 1, given in some cells indicated in red.}
\label{fig:mesh-obs}
\end{figure}

\begin{figure}[hbt!]
\centering
\begin{tabular}{cc}
\includegraphics[width=0.4 \textwidth]{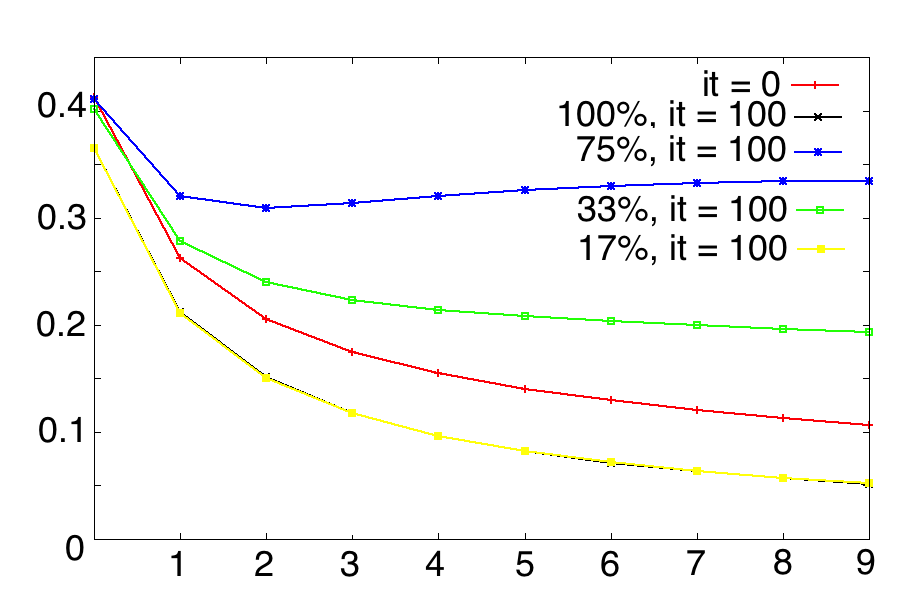}&\includegraphics[width=0.4 \textwidth]{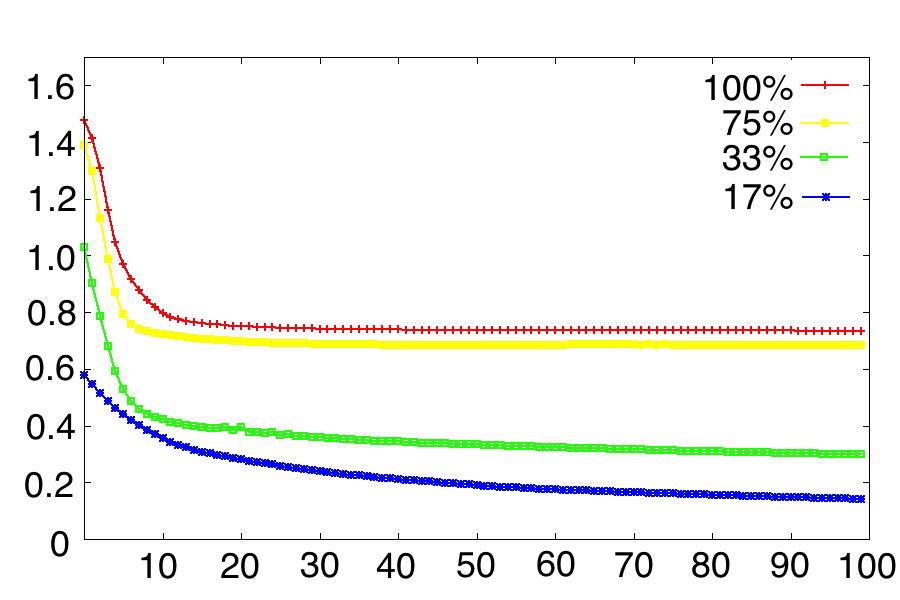}\\
(a) $L^2$ error v.s. time & (b) $J$ v.s. iteration index
\end{tabular}
\caption{Numerical results for Case 1 using 4 multiscale basis functions, with different amount of observational
data $g^K$ in some cells shown in Fig. \ref{fig:mesh-obs}:  (a) the $L^2$ error for cell average v.s. time $t$,
and (b) the functional value $J$ v.s. the iteration index.}
\label{fig:inv3r-obs}
\end{figure}

Next we illustrate the sensitivity of the numerical results by the multiscale inversion algorithm with respect to various
algorithmic parameters. In Fig. \ref{fig:inv3-ad} we present the result for two different step lengths $\epsilon=10^{-12}$
and $\epsilon=10^{-13}$, where the mass and stiffness matrices are updated adaptively in selected regions and also over
the whole computational domain. One observes that the errors and residuals are comparable when the iteration reaches
convergence, but with a larger step size can greatly speedup the convergence of the algorithm (whenever it does not violate the
step size restriction, as usual for gradient descent type algorithms). Further, the results for the adaptive local update and
all cells update of the mass matrices are comparable with each other. Thus the inversion with only local update in the
selected regions affects little the reconstruction results. However, numerically, we observe that the local update is much
more stable than the global update, e.g., a larger step size $\epsilon$, due to the fact that the local update involves much
few unknowns. In Fig. \ref{fig:inv3r-rand}, we present the numerical results for the case of observational data contaminated
with different amount of noise. It is observed that the results are fairly stable with respect to the present of data noise,
up to 100 iterations, due to the priors we specified on the discrete parameters, clearly indicating the stability of the
regularized formulation. Naturally, the error and residual increases with the noise level. Last, our inversion algorithm
essentially employs local data to update the local coarse grid directly, and thus it is expected that the algorithm can work
well as long as the related local data over the interested region is available. This is confirmed by the numerical results in
Figs. \ref{fig:mesh-obs} and \ref{fig:inv3r-obs}.  However, with sparser data available, a
stronger regularization is needed to maintain the stability of the algorithm, and more informative priors, e.g., sparsity
or total variation, may be imposed \cite{JinMaass:2012,SchusterKaltenbacher:2012,ItoJin:2015}.

%

\section{Conclusions}

In this work, we have developed a generalized multiscale inversion algorithm for heterogeneous problems. It is based
on the generalized multiscale finite element method (GMsFEM), where one constructs multiscale basis functions to
capture the non-localizable features, and the algorithm assumes that the problem admits a reduced-order model on a
coarse grid. Then, instead of seeking coarse-grid permeabilities, we seek the discretization parameters that are obtained
from the GMsFEM formulation. Our approaches are especially suitable for problems with fractures or high-conductivity channels,
when upscaling the permeability can result in very large errors. Thus, it is important to consider a more general
multiscale approach. In our approach, we do not compute multiscale basis functions and do not recover the fine-scale
permeability field. Instead, we compute the averaged coarse-grid discretization parameters, i.e., integrated responses corresponding to
unknown multiscale basis functions. We have discussed various regularizations and a Bayesian framework, as well
as the important ingredients of the algorithm, and illustrated the approach with numerical results for fractured media.
Our numerical experiments clearly illustrate the feasibility and significant potential of the approach for inverse
problems for heterogeneous problems, and it motivates a rigorous analysis of the proposed approach.

\bibliographystyle{abbrv}
\bibliography{msinv}
\end{document}